\documentclass[letterpaper, 10 pt, conference]{ieeeconf}
\IEEEoverridecommandlockouts
\usepackage{cite}
\usepackage{amsmath,amssymb,amsfonts}
\usepackage{algorithmic}
\usepackage{graphicx}
\usepackage{textcomp}
\usepackage{xcolor}
\usepackage{verbatim}
\usepackage{subcaption}
\usepackage{epstopdf}
\def\BibTeX{{\rm B\kern-.05em{\sc i\kern-.025em b}\kern-.08em
    T\kern-.1667em\lower.7ex\hbox{E}\kern-.125emX}}

\title{\LARGE \bf Thermal management optimization of Battery Electric Vehicles via hierarchical NMPC with CMO-trained Neural Models
}
\author{
Francesco Ripa\textsuperscript{a}, Diego Regruto\textsuperscript{a}, Pasquale Ceres\textsuperscript{b}, Christopher Strano\textsuperscript{b} \thanks{\textsuperscript{a} F. Ripa and D. Regruto are with Department of Control and Computer Engineering, Politecnico di Torino, Italy;  \textsuperscript{b} P. Ceres and C. Strano are with Fiat Research Center (CRF), Stellantis N.V., Turin, Italy}}
\begin{document}
\maketitle
\thispagestyle{empty}
\pagestyle{empty}

\begin{abstract}
This paper addresses the problem of energy-efficient thermal management in Battery Electric Vehicles (BEVs), with the goal of extending driving range while ensuring passenger comfort and battery thermal safety. A hierarchical Nonlinear Model Predictive Control (NMPC) framework is proposed, in which the supervisory layer explicitly balances energy consumption and cabin temperature reference tracking over a long horizon and translates this trade-off into optimal actuator command sequences for the integrated heating, ventilation, and air conditioning (HVAC) system and the battery thermal management (BTM) system. These command sequences are then enforced by a lower-layer NMPC, which introduces local corrections to compensate disturbances and modeling uncertainty while tracking an externally specified cabin temperature reference trajectory. To enable predictive control of the highly nonlinear and strongly coupled thermal dynamics, data-driven neural network models are employed as prediction models within both control layers. The proposed approach is developed and validated using a high-fidelity BEV thermal simulator developed at Politecnico di Torino in collaboration with Stellantis N.V., and benchmarked against an industrial rule-based (RB) control strategy over the WLTC driving cycle. Simulation results demonstrate a significant reduction in energy consumption compared to the baseline strategy, while satisfying comfort and battery temperature constraints.
\end{abstract}

%\begin{IEEEkeywords}
%Energy saving, model predictive control, Machine Learning 
%\end{IEEEkeywords}

\section{Introduction}

In the last decade, increasing awareness of climate change, air 
pollution, and the progressive depletion of fossil fuels has 
significantly accelerated the transition toward electric vehicles 
(EVs). This shift is largely driven by the promise of reduced 
environmental impact, together with lower maintenance and operating 
costs compared to conventional internal combustion engine vehicles 
\cite{kumar2020adoption,choi2018effect}. Despite the progress, 
driving range remains a key barrier, strongly affected by ambient 
temperature through the combined impact on battery performance, 
powertrain efficiency, and cabin climate control demand. Recent 
chassis-dynamometer tests show that at $-18^\circ$C, usable battery 
energy decreases by 4--8\% compared to $22^\circ$C and the overall 
driving range can drop by up to 60\%, with HVAC energy consumption 
being one of the most relevant contributing factors 
\cite{seo2025effects}. These results strongly motivate the need for 
energy-efficient thermal management strategies, particularly under 
winter and summer operating conditions.

Thermal management in EVs serves a dual purpose: ensuring battery 
performance, safety, and longevity, while maintaining acceptable 
passenger comfort through cabin thermal regulation. State-of-practice 
solutions rely on rule-based (RB) supervisory control logic, where 
cooling and heating actions are triggered based on predefined 
temperature thresholds \cite{basu2019battery,kummer2013method}. 
Beyond rule-based solutions, optimization-based methods have been 
extensively investigated. MPC-based approaches have been proposed for 
cabin HVAC control \cite{wang2019mpc} and battery thermal management 
\cite{bauer2014thermal,masoudi2015battery,masoudi2017mpc,
tao2014cooling,tao2015hybrid}, and joint cabin--battery thermal 
management has been studied in hierarchical MPC settings 
\cite{amini2019cabin}. Despite these advances, predictive control of 
integrated systems remains challenging when actions are applied 
directly at the actuator level, as strong coupling among heat pump 
dynamics, refrigerant circuits, and routing valves often necessitates 
simplifications that limit modeling accuracy over wide operating 
conditions. Large thermal inertia and slow dynamics further require 
long prediction horizons, while predictive information introduces uncertainty that must be balanced against real-time 
constraints. These challenges naturally motivate hierarchical control 
architectures separating long-term planning from short-term regulation 
\cite{scattolini2009architectures,amini2020hierarchical}.

In this work, we propose a two-level hierarchical NMPC strategy for 
the integrated management of cabin HVAC and battery thermal systems 
under both summer and winter operating conditions, with three main 
contributions relative to the existing literature. First, unlike 
conventional hierarchical MPC schemes that generate thermal state 
references at the supervisory level, the proposed upper layer 
explicitly optimizes actuator command trajectories over a long horizon 
exploiting speed profile prediction, while the lower-layer NMPC tracks these 
trajectories with local corrections for disturbance rejection and 
cabin temperature regulation. This input-scheduling paradigm enables 
tighter coordination between long-horizon energy planning and 
real-time control. Second, compact recurrent neural network (RNN) 
models trained via the Controller Multiplier Optimization (CMO) 
approach of \cite{cerone2025new} are consistently embedded in both 
control layers as prediction models: the CMO training mitigates 
vanishing and exploding gradient issues and yields accurate 
control-oriented models without requiring access to the internal 
equations of the high-fidelity GT-SUITE simulator. Third, the 
proposed strategy is validated against a Stellantis rule-based 
industrial benchmark over the WLTC driving cycle, demonstrating 
energy savings exceeding 12\% under both cold and hot ambient 
conditions while satisfying all thermal and actuator constraints.

The remainder of the paper is organized as follows: 
Section~\ref{sec:sys_descr} describes the system and its modeling; 
Section~\ref{sec:math_form} presents the hierarchical NMPC framework 
and the neural network identification procedure; 
Section~\ref{sec:simulation} reports simulation results; 
Section~\ref{sec:conclusion} draws conclusions.

\section{System description} \label{sec:sys_descr}

The control strategy developed in this work targets the thermal 
system architecture defined within the OPTEMUS framework. In the 
context of this paper, OPTEMUS refers to a specific integrated 
thermal management architecture conceived to simultaneously regulate 
cabin air conditioning and battery temperature, exploiting thermal 
synergies among vehicle subsystems to improve energy efficiency.

At the core of the thermal management system lies the Compact 
Refrigeration Unit (CRU), a heat-pump-based unit responsible for 
generating and distributing heating or cooling capacity according to 
the thermal requirements of the vehicle components. When operating 
conditions permit, the CRU transfers heat extracted from electrified 
powertrain components, including the battery, to support cabin 
heating. In scenarios where the cabin and battery simultaneously 
require heating or cooling, the CRU activates an appropriate 
combination of heat sources and sinks to deliver the required thermal 
power. Refrigerant routing is managed through multiple three-way 
valves, which distribute hot and cold refrigerant flows between the 
cabin HVAC circuit and the battery thermal loop depending on 
instantaneous operating conditions.

\begin{figure}[t]
  \centering
  \includegraphics[width=0.5\textwidth]{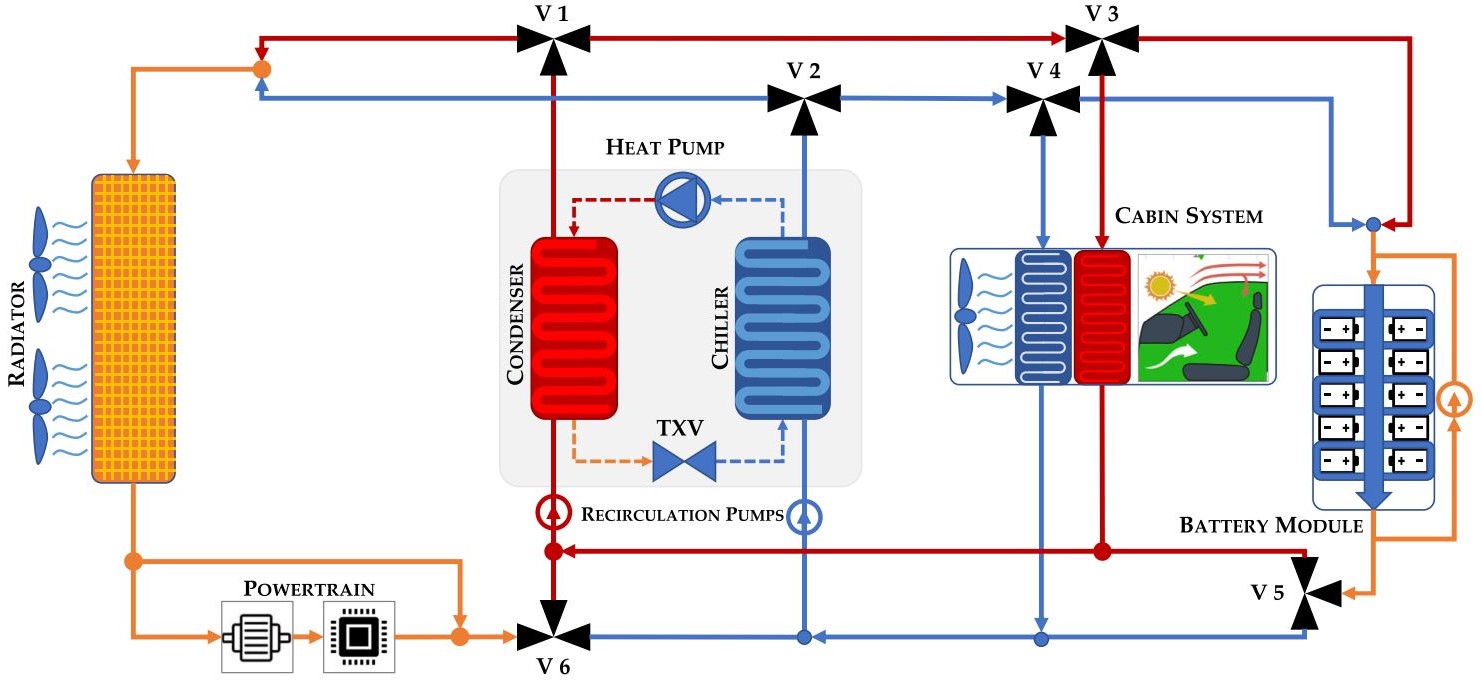}
  \caption{Thermal management system layout considered in this work 
  \cite{rolando2024virtual}.}
  \label{fig:layout}
\end{figure}

A detailed schematic of the thermal management system is shown in 
Fig.~\ref{fig:layout}, developed in the GT-SUITE software environment 
and presented in \cite{rolando2024virtual}. The coupling between 
cabin and battery thermal dynamics motivates a coordinated control 
strategy accounting for energy efficiency and operational constraints. 
The system behavior is described, in discrete time, by the state 
vector
\begin{equation}
x(k)=\big[T_{\mathrm{bat}}(k)\ \ SOC(k)\ \ \Delta T(k)\ \ 
T_{\mathrm{cab}}(k)\big]^{\top},
\label{eq:state_vector}
\end{equation}
where $T_{\mathrm{bat}}$ is the battery temperature, $SOC$ is the 
battery state of charge, $T_{\mathrm{cab}}$ is the cabin temperature, 
and $\Delta T$ is the battery temperature differential, defined as 
the thermal excursion between coolant inlet temperature and battery 
temperature, included among the constrained states as persistently 
large values indicate limited thermal regulation effectiveness and 
may accelerate battery degradation. Under the assumption of full 
state availability, the measured output is
\begin{equation}
y(k)=g(x(k))=x(k).
\label{eq:output_vector}
\end{equation}
All state variables are considered measurable in practice: 
$T_{\mathrm{bat}}$ and $SOC$ are standard battery management system 
(BMS) outputs, $T_{\mathrm{cab}}$ is measured by the climate control 
unit, and $\Delta T$ is derived from coolant inlet and battery 
temperature sensors. The development of a state observer, aimed at reducing the number of required sensors, is left for future investigations.

The CRU controls the thermal response through compressor actuation 
and refrigerant routing. The proposed controller directly regulates 
low-level actuators, namely the compressor speed and a main routing 
valve, rather than high-level thermal power requests. Accordingly, 
the manipulated input vector is %\vspace{-0.1cm}

\begin{equation}
u(k)=\big[v_{\mathrm{comp}}(k)\ \ \%\mathrm{valve}(k)\big]^{\top},
\label{eq:input_vector}
%\vspace{-0.1cm}
\end{equation}
where $v_{\mathrm{comp}}$ is the compressor speed and 
$\%\mathrm{valve}$ denotes the opening ratio of the three-way valve 
in $[0,1]$, where 1 directs flow entirely toward the cabin HVAC 
circuit, 0 toward the battery thermal management circuit, and 
intermediate values split the flow proportionally. The controller 
commands $V_3$ in cooling mode and $V_4$ in heating mode; the 
remaining valves ($V_1$--$V_2$--$V_5$--$V_6$) and the cabin fan are 
managed by internal low-level logic that configures the refrigerant 
path according to the selected operating mode.

Two distinct prediction models are derived for cooling and heating 
operation respectively. A data-driven approach based on recurrent 
neural networks is adopted rather than a first-principles model, 
as the complexity of the refrigerant circuit dynamics and the 
proprietary nature of the GT-SUITE simulator make physics-based 
modeling impractical at the required level of detail. RNN models 
allow compact, control-oriented predictors to be extracted directly 
from simulation data, capturing the strongly nonlinear and coupled 
thermal behavior across both operating regimes, as described in 
Section \ref{RNN}.

\section{Problem Statement and Proposed Control Strategies}
\label{sec:math_form}

The integrated regulation of the cabin system and the battery cooling and heating circuit gives rise to a multi-objective control problem characterized by competing performance requirements and strong internal couplings. Passenger comfort requires the cabin temperature $T_{\mathrm{cab}}$ to accurately track a prescribed reference $T_{\mathrm{cab,ref}}$, which is generated based on the ambient temperature, the user-selected setpoint, and comfort-oriented constraints on the rate of temperature variation. More specifically,  $T_{\mathrm{cab,ref}}(k)$ is generated by filtering the step between the initial cabin temperature and the desired setpoint, so as to obtain a smooth profile with prescribed transient characteristics (e.g., rise/settling time) consistent with passenger comfort requirements.
At the same time, battery thermal safety and operability require that the battery temperature $T_{\mathrm{bat}}$ and the associated thermal excursion $\Delta T$ remain within admissible bounds. These objectives must be pursued while minimizing the energy consumption of the thermal management system, which directly impacts the driving range of the vehicle.

%With the state vector $x(k)$ and the input vector $u(k)$ defined in \eqref{eq:state_vector}--\eqref{eq:input_vector}, the system evolution is described in discrete time by a nonlinear state-space model of the ,
%\begin{equation}
%x(k+1)=f_{\mathrm{NN}}\!\big(x(k),u(k)\big),
%\label{eq:dynamics_ps}
%\end{equation}
%which captures the nonlinear and coupled thermal dynamics of the system. 

The overall performance index to be optimized is defined as
\begin{equation}
J_s = J_{1}+J_{2},
\label{eq:J_sum}
\end{equation}
where
\begin{equation}
J_{1}=\sum_{i=0}^{p-1}\left\|T_{\mathrm{cab}}(k{+}i|k)-T_{\mathrm{cab,ref}}(k+i)\right\|_{Q}^{2}
\label{eq:J1}
\end{equation}
penalizes deviations from the desired cabin-temperature reference trajectory, and the term $J_2$ captures the energy consumption of the thermal management system. Since the thermal power demand can be expressed through the absorbed current $i_{\mathrm{abs}}$ which, in turn, depends on the compressor speed $v_{\mathrm{comp}}$ and the valve activation, the energy consumption objective $J_2$ is defined as the weighted energy of the command input defined in equation (3), i.e.:
\begin{equation}
J_{2}=\sum_{i=0}^{p-1}\left\|u(k{+}i|k)\right\|_{R}^{2}.
\label{eq:J2}
\end{equation}
The weighted quadratic norms in equation (5) and (6) are defined as $\|z\|_{Q}^{2}:=Qz^{2}$ for scalar $z$ and $\|w\|_{R}^{2}:=w^{\top}Rw$ for vector $w$, with $Q>0$ and $R\succ 0$.
From both an optimization and an implementation perspective, the solution of this control problem is challenging. The presence of large thermal inertia and slow temperature dynamics requires long prediction horizons in order to appropriately shape the thermal response and exploit energy-efficient operating conditions. Moreover, traction and thermal loads are intrinsically coupled, as the battery represents the sole energy source for both propulsion and auxiliary systems, and variations in thermal demand directly affect the available energy for driving. Additional complexity arises from the nonlinearities introduced by the refrigerant routing logic of the compact refrigeration unit, which couples valve configurations with power demand and heat-flow distribution. These considerations motivate the adoption of a hierarchical MPC architecture, in which long-horizon planning and short-horizon regulation are explicitly separated. In this framework, preview information such as speed profiles is exploited to compute energy-efficient thermal actions over an extended horizon, while a lower-level controller ensures real-time feasibility and constraint satisfaction.

\subsection{Proposed hierarchical NMPC strategy}
\label{subsec:hnmpc}

The proposed control strategy is structured as a two-layer hierarchical nonlinear model predictive control (NMPC) scheme, designed to address the different time scales and objectives characterizing integrated cabin HVAC and battery thermal management. The two control layers operate at different sampling rates and fulfill complementary roles within the overall architecture.

The supervisory layer is responsible for long-horizon planning of the thermal management actions. At each supervisory update, a finite-horizon optimal control problem is solved using preview information provided by the navigation system (and/or vehicle-to-vehicle (V2V) and vehicle-to-infrastructure (V2I) sources), in particular the predicted vehicle speed profile, which affects the thermal dynamics. The supervisory optimization explicitly balances energy consumption and cabin temperature tracking over a long prediction horizon, while accounting for battery thermal constraints. The outcome of this optimization is a scheduled actuation policy expressed in terms of optimal command input trajectories for the thermal actuators. The supervisory problem is solved at a slower time scale and is periodically updated in a receding-horizon fashion, generating updated command trajectories that are supplied as references to the lower control layer.

The lower-layer NMPC operates at a faster time scale and is tasked with enforcing the scheduled actuation policy. Specifically, the lower layer tracks the command input trajectories provided by the supervisory layer while introducing local corrections to compensate for disturbances, modeling uncertainty, and unmodeled dynamics. In addition, an output-tracking term is included in the lower-layer cost function to ensure tracking of an externally specified cabin temperature reference trajectory. In this way, the lower layer does not re-optimize the long-horizon energy/comfort trade-off determined by the supervisory layer, but rather enforces the planned actuation policy and locally corrects it to maintain the desired thermal behavior.

Both control layers employ the same neural-network-based model (discussed in Section III.B) to describe the complex nonlinear system dynamics and to perform state prediction, ensuring consistency between long-horizon planning and short-horizon regulation. The high-fidelity GT-SUITE vehicle thermal model is used exclusively for closed-loop simulation and validation of the proposed control strategy and is not embedded within the predictive controllers. \\

\subsubsection{Supervisory layer (long-horizon optimization with preview)}

Let $p_{\mathrm{s}}$ denote the scheduling horizon and define the stacked decision vector
\begin{equation}
U_{\mathrm{s}}(k)=
\begin{bmatrix}
u_{\mathrm{s}}(k|k)^\top\!\! &
u_{\mathrm{s}}(k{+}1|k)^\top\!\! &
\cdots\!\! &
u_{\mathrm{s}}(k{+}p_{\mathrm{s}}{-}1|k)^\top
\end{bmatrix}^{\!\top}
\label{eq:Us_def}
\end{equation}
The optimal control problem is formulated as
\begin{IEEEeqnarray}{rCl}
U_{\mathrm{s}}^{*}(k) &=&
\arg\min_{U_{\mathrm{s}}(k)}\ J_{\mathrm{s}}\!\left(x(k),U_{\mathrm{s}}(k)\right)
\label{eq:opt_sched}\\
\text{s.t.}\quad
x(k{+}i{+}1|k) &=&
f_{\mathrm{NN}}\!\left(x(k{+}i|k),u_{\mathrm{s}}(k{+}i|k),v(k{+}i|k)\right),
\nonumber\\
&& i=0,\ldots,p_{\mathrm{s}}{-}1,
\label{eq:opt_sched_dyn}\\
x_{\min} \leq x(k{+}i|k) &\leq& x_{\max},
\quad i=0,\ldots,p_{\mathrm{s}},
\label{eq:opt_sched_xc}\\
u_{\min} \leq u_{\mathrm{s}}(k{+}i|k) &\leq& u_{\max},
\quad i=0,\ldots,p_{\mathrm{s}}{-}1.
\label{eq:opt_sched_uc}
\end{IEEEeqnarray}
where $f_{\mathrm{NN}}$ in equation (9) is the neural network model describing the system dynamics.
The cost $J_{\mathrm{s}}$ encodes the long-horizon energy objective and the cabin comfort requirement. The speed profile preview $v(k{+}i|k)$, provided by the navigation system, is assumed a known exogenous input to the function $f_{\mathrm{NN}}$. The resulting optimal sequence defines the scheduled actuator trajectory provided to the lower piloting layer. \\

\subsubsection{Piloting layer (short-horizon tracking NMPC)}

Let $p_{\mathrm{p}}$ denote the piloting horizon, with $p_{\mathrm{p}}\ll p_{\mathrm{s}}$, and define
\begin{equation}
U_{\mathrm{p}}(k)=\left\{u_{\mathrm{p}}(k|k),u_{\mathrm{p}}(k+1|k),\dots,u_{\mathrm{p}}(k+p_{\mathrm{p}}-1|k)\right\}.
\label{eq:Up_def}
\end{equation}

\begin{comment}
    The real-time tracking problem solved by the piloting layer is formulated as
\begin{equation}
\begin{aligned}
U_{\mathrm{p}}^{*}(k)=\ &\arg\min_{U_{\mathrm{p}}(k)}\ J_{\mathrm{p}}\!\left(x(k),U_{\mathrm{p}}(k)\right)\\
\text{s.t.}\ \ &x_{k+i+1|k}=f_{\mathrm{NN}}\!\left(x_{k+i|k},u_{\mathrm{p},k+i|k}\right),\\
&\hspace{2.45cm} i=0,\dots,p_{\mathrm{p}}-1,\\
&x_{\min}\le x_{k+i|k}\le x_{\max},\ \ i=0,\dots,p_{\mathrm{p}},\\
&u_{\min}\le u_{\mathrm{p},k+i|k}\le u_{\max},\ \ i=0,\dots,p_{\mathrm{p}}-1,
\end{aligned}
\label{eq:opt_pilot}
\end{equation}
where the cost function is given by
\begin{equation}
\begin{aligned}
J_{\mathrm{p}}=&\ \left\|T_{\mathrm{cab}}^{p_{\mathrm{p}}}-T_{\mathrm{ref}}^{p_{\mathrm{p}}}\right\|_{P_f}^{2}
+\sum_{i=0}^{p_{\mathrm{p}}-1}\Big(
\left\|T_{\mathrm{cab}}^{i}-T_{\mathrm{ref}}^{i}\right\|_{Q}^{2}\\
&\hspace{2.6cm}+
\left\|u_{\mathrm{p}}^{i}-u_{\mathrm{s}}^{i}\right\|_{R}^{2}
\Big),
\end{aligned}
\label{eq:pilot_cost}
\end{equation}
with $T_{\mathrm{cab}}^{i}:=T_{\mathrm{cab}}(k{+}i|k)$, $T_{\mathrm{ref}}^{i}:=T_{\mathrm{cab,ref}}(k{+}i|k)$, $u_{\mathrm{p}}^{i}:=u_{\mathrm{p}}(k{+}i|k)$, and $u_{\mathrm{s}}^{i}:=u_{\mathrm{s}}^{*}(k{+}i|k)$. {\color{orange} The notation $\|z\|_{M}^{2}=z^{\top}Mz$ is adopted, where $Q\succeq 0$, $R\succ 0$, and $P_f\succeq 0$ are weighting matrices.} \\
\end{comment}

The optimal control problem solved by the piloting layer is
\begin{IEEEeqnarray}{rCl}
U_{\mathrm{p}}^{*}(k) &=& \arg\min_{U_{\mathrm{p}}(k)}\ J_{\mathrm{p}}\!\left(x(k),U_{\mathrm{p}}(k)\right)
\label{eq:opt_pilot}\\
\text{s.t.}\quad x(k{+}i{+}1|k) &=& f_{\mathrm{NN}}\!\left(x(k{+}i|k),u_{\mathrm{p}}(k{+}i|k),v(k{+}i|k)\right),
\nonumber\\
&& i=0,\dots,p_{\mathrm{p}}-1, \label{eq:opt_pilot_dyn}\\
x_{\min} \le x(k{+}i|k) &\le& x_{\max}, \quad i=0,\dots,p_{\mathrm{p}}, \label{eq:opt_pilot_x}\\
u_{\min} \le u_{\mathrm{p}}(k{+}i|k) &\le& u_{\max}, \quad i=0,\dots,p_{\mathrm{p}}-1.
\label{eq:opt_pilot_u}
\end{IEEEeqnarray}
where the cost function to be minimized
\begin{equation}
J_{\mathrm{p}}=\sum_{i=0}^{p_{\mathrm{p}}-1}\Big(\|e_T(k{+}i|k)\|_{Q}^{2}+\|e_u(k{+}i|k)\|_{R}^{2}\Big),
\label{eq:pilot_cost}
\end{equation}
is a weighted combination of the tracking errors defined as follows:
\begin{equation}
e_T(k{+}i|k)\triangleq T_{\mathrm{cab}}(k{+}i|k)-T_{\mathrm{cab,ref}}(k+i),
\label{eq:eT}
\end{equation}
\begin{equation}
e_u(k{+}i|k)\triangleq u_{\mathrm{p}}(k{+}i|k)-u_{\mathrm{s}}^{*}(k{+}i|k).
\label{eq:eu}
\end{equation}

The weighted quadratic norm in (17) are defined as  $\|z\|_{Q}^{2}:=Qz^{2}$ for scalar $z$ and $\|w\|_{R}^{2}:=w^{\top}Rw$ for vector $w$.
Here, $Q>0$ weights cabin-temperature tracking and $R\succ 0$ weights the deviation from the command inputs trajectories scheduled by the supervisory layer.
Also in the piloting layer the speed preview $v(k{+}i|k)$ is a known exogenous input to $f_{\mathrm{NN}}$.

%\subsubsection{Piloting layer (short-horizon NMPC)}
The piloting layer operates online with a short prediction horizon $p_p \ll p_s$ and is responsible for real-time regulation, constraint enforcement, and disturbance rejection. At each sampling instant $k$, the controller solves the NMPC problem in (15) and (16), where the scheduled input trajectory $u_s^{*}(\cdot|k)$ generated by the supervisory layer acts as a reference signal for the optimization. Differently from standard hierarchical MPC formulations based on output/states reference tracking, the hierarchical reference enters the piloting-layer problem through the quadratic input-tracking term $\|u_p - u_s\|_R^2$, which provides the primary grounding of the optimization to the scheduled actuator trajectory. Output variables are still included in the cost to ensure cabin temperature regulation and constraint satisfaction, but they do not act as the main driver of the optimization. In the piloting layer the two cost terms play distinct roles: the input-tracking term $\|e_u\|_Q^2$ defines the nominal actuation policy, while the output-tracking term $\|e_T\|_R^2$ introduces the local input corrections required to track the externally specified cabin temperature reference.

At each sampling instant $k$, the NMPC problem returns an optimal input sequence $U_p^{\*}(k)=\{u_p^{\*}(k|k),\ldots,u_p^{\*}(k+p_p-1|k)\}$; only the first element is applied to the plant and the optimization is repeated at the next sampling instant, according to standard receding-horizon strategy.

The scheduled input trajectory computed by the upper layer thus acts as a reference signal for the online tracking problem, ensuring consistency between long-term planning and short-term regulation. At each sampling instant $k$, the NMPC problem returns an optimal input sequence
$U_{\mathrm{p}}^{*}(k)=\{u_{\mathrm{p}}^{*}(k|k),\dots,u_{\mathrm{p}}^{*}(k{+}p_{\mathrm{p}}{-}1|k)\}$; only the first element is applied to the plant and the optimization is repeated at $k{+}1$.

% Subsection B will describe how f_NN is constructed and trained.
\subsection{Recurrent Neural Network based modeling of the thermal dynamics}\label{RNN}
A core element of the proposed hierarchical NMPC framework is the prediction model, which must capture the nonlinear and coupled thermal dynamics of the integrated system. In this work, the nonlinear dynamical model is learned via black-box system identification using recurrent neural networks (RNNs), trained on simulation data obtained by running the high-fidelity GT-SUIT model. Two separate models are identified to account for the different thermal behaviors observed in cooling (hot summer) and heating (cold winter) operating conditions, yielding $f_{\mathrm{NN,HOT}}$ and $f_{\mathrm{NN,COLD}}$, respectively.

The identification problem is formulated by applying the approach proposed in Sec.~6.2 of \cite{cerone2025constrained}. More precisely, the following state-space structure is considered
\begin{equation}
x_{t+1}=f_{NN}(x_t,u_t,v_t,\theta),\qquad
y_t=x_t,
\label{eq:ss_cmo}
\end{equation}
where $f_{NN}$ is a recurrent neural network of suitable structure (see \cite{cerone2025constrained} for details),  $\theta$ is the vector of model parameters to be trained, and $v_t$ denotes the (measured/previewed) vehicle speed provided by the navigation system and  treated as a known exogenous input to the nonlinear dynamics. Following the approach proposed in \cite{cerone2025constrained}, the problem of estimating the parameter $\theta$ can be recast as the constrained optimization problem 
\vspace{-0.5cm}
\begin{equation}
\begin{aligned}
\min_{\theta,x_1,\dots,x_N}\ \ &\sum_{t=1}^{N}\big\|x_t-\tilde{y}_t\big\|_{2}^{2}+\rho(\theta)\\
\text{s.t.}\ \ &x_{t+1}=f_{NN}(x_t,u_t,v_t,\theta),\qquad t=1,\dots,N-1,
\end{aligned}
\label{eq:cmo_problem}
\end{equation}
which can be efficiently solved using CMO approach proposed in \cite{cerone2025new} which overcome common issues affecting the training of recurrent networks (vanishing/exploding gradient problems). We point the interested reader to the work \cite{cerone2025constrained} for a detailed discussion. 

Two RNN models $f_{NN,HOT}$ and $f_{NN,COLD}$  are trained on data obtained by simulating the high-fidelity GT-SUIT model under cold winter (heating) and hot summer (cooling) conditions respectively.The selection between the two models is determined by the 
current operating mode of the thermal system: $f_{\mathrm{NN,COLD}}$ 
is active during heating operation and $f_{\mathrm{NN,HOT}}$ during 
cooling operation, consistently with the low-level valve logic that 
configures the refrigerant circuit. Both predictors adopt compact architectures to preserve control-oriented simplicity: $f_{\mathrm{NN,COLD}}$ employs one recurrent layer with $4$ neurons, while $f_{\mathrm{NN,HOT}}$ employs one recurrent layer with $6$ neurons. Model performance is assessed through the Root Mean Square Error (RMSE),
\begin{equation}
\mathrm{RMSE}=\sqrt{\frac{1}{n}\sum_{j=1}^{n}\big(y_j-\hat{y}_j\big)^2}.
\label{eq:rmse}
\end{equation}

\begin{table}[t]
\centering
\caption{RMSE of the identified RNN predictors in normalized domain.}
\label{tab:rmse_summary}
\renewcommand{\arraystretch}{1.15}
\begin{tabular}{lcccc}
\hline
 & \multicolumn{2}{c}{\textbf{Cold case}} & \multicolumn{2}{c}{\textbf{Hot case}}\\
\textbf{State} & \textbf{Train} & \textbf{Val} & \textbf{Train} & \textbf{Val}\\
\hline
$T_{\mathrm{bat}}$ & 0.043 & 0.100 & 0.003 & 0.019 \\
$T_{\mathrm{cab}}$ & 0.110 & 0.170 & 0.042 & 0.140 \\
$\Delta T$         & 0.067 & 0.072 & 0.032 & 0.091 \\
$SOC$              & 0.014 & 0.018 & 0.008 & 0.044 \\
\hline
\end{tabular}
\end{table}

Table~\ref{tab:rmse_summary} reports the RMSE values obtained for each state component in both operating modes. The \emph{Train} columns refer to the dataset used for parameter identification, whereas the \emph{Val} columns are computed on an independent validation dataset not used in the identification step. The relatively small gap between training and validation errors support the use of these compact models in the NMPC scheme as sufficiently accurate prediction models.

\section{Simulation Results}
\label{sec:simulation}

The proposed hierarchical NMPC strategy was validated in MATLAB/Simulink using the high-fidelity thermal plant model developed at the Politecnico di Torino in collaboration with Stellantin N.V. (presented in \cite{rolando2024virtual}), and benchmarked against the Stellantis rule-based (RB) controller. The RNN prediction models described in Sec.~\ref{RNN} are used as control-oriented models within the hierarchical NMPC optimization scheme.

\begin{figure}[t]
  \centering
  \includegraphics[width=0.9\columnwidth]{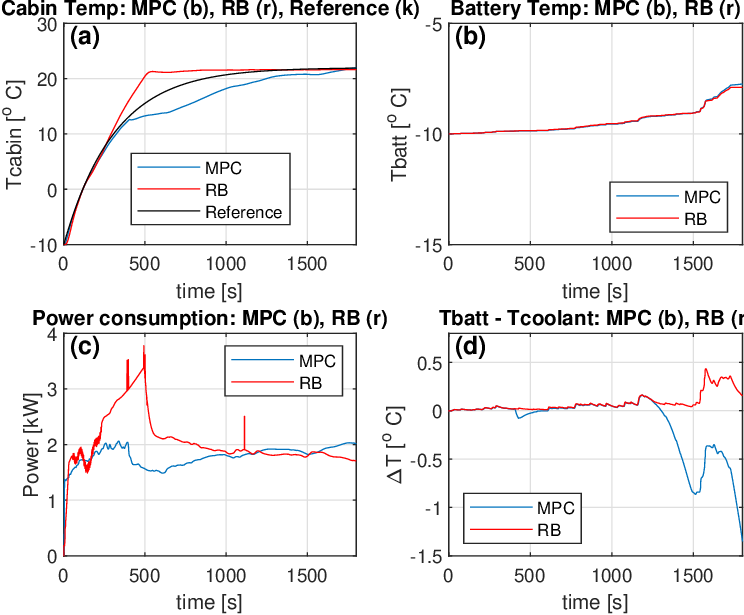}
  \caption{Cold case: state comparison between MPC (blue) and RB (red): $T_{\mathrm{cab}}$, $T_{\mathrm{bat}}$, $Power$ $consumption$, and $\Delta T$.}
  \label{fig:cold_states}
\end{figure}

\begin{figure}[t]
  \centering
  \includegraphics[width=\columnwidth]{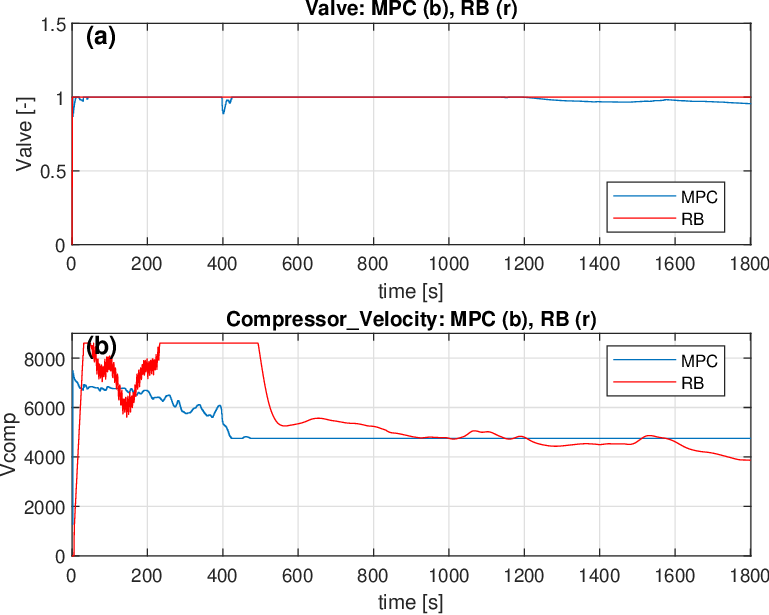}
  \caption{Cold case: input comparison between MPC (blue) and RB (red): valve command and compressor speed.}
  \label{fig:cold_inputs}
\end{figure}

\begin{figure}[t]
  \centering
  \includegraphics[width=\columnwidth]{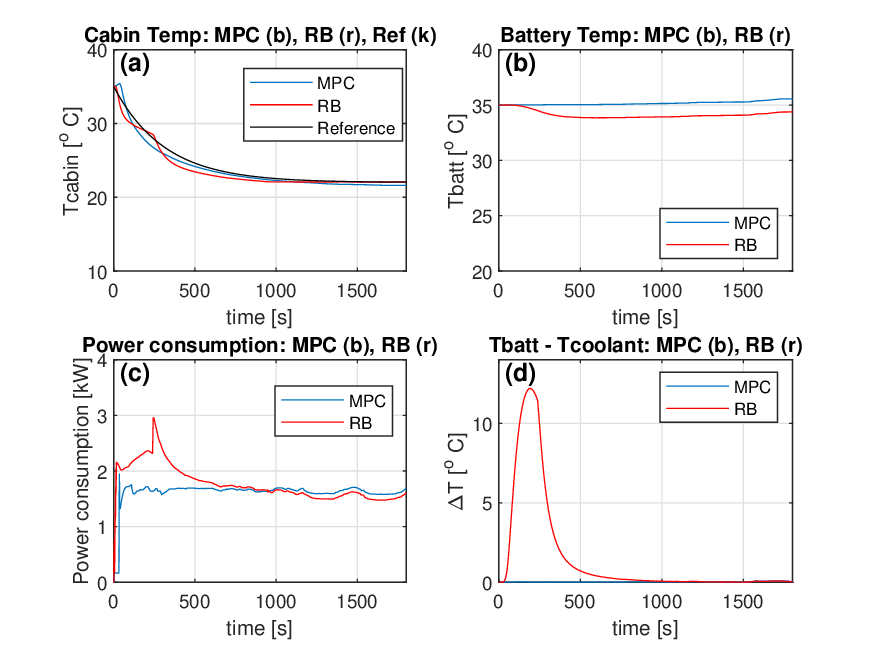}
  \caption{Hot case: state comparison between MPC (blue) and RB (red): $T_{\mathrm{cab}}$, $T_{\mathrm{bat}}$, $Power$ $consumption$, and $\Delta T$.}
  \label{fig:hot_states}
  \vspace{-0.5cm}
\end{figure}

\begin{figure}[t]
  \centering
  \includegraphics[width=\columnwidth]{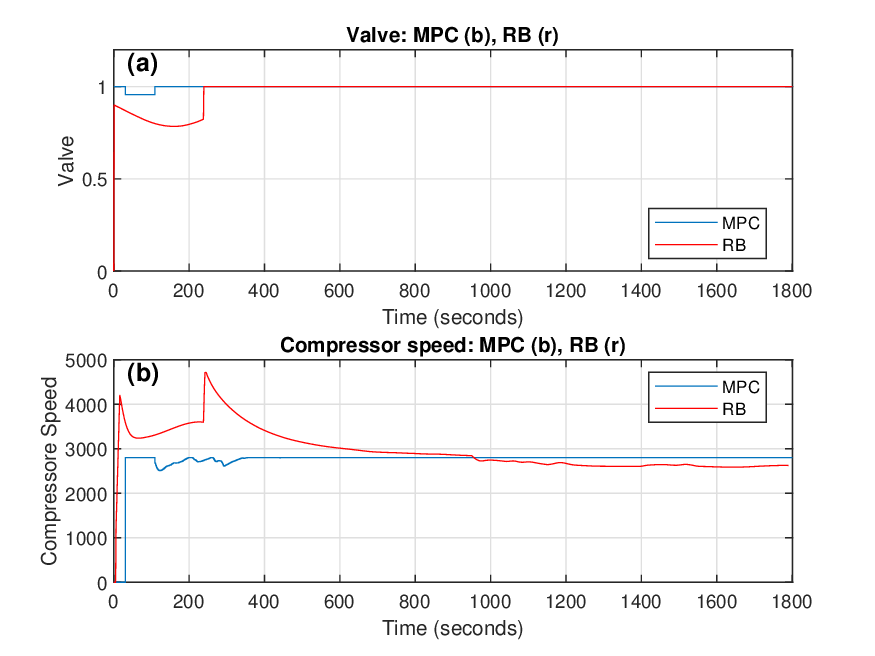}
  \caption{Hot case: input comparison between MPC (blue) and RB (red): valve command and compressor speed.}
  \label{fig:hot_inputs}
  \vspace{-0.5cm}
\end{figure}

Two representative operating conditions are considered. In the \textit{cold winter case}, the external ambient temperature is set to $-10^{\circ}$C and the objective is to heat the cabin to the final target set point of $22^{\circ}$C by tracking the reference trajectory $T_{\mathrm{cab,ref}}$ depicted as a black line in Fig. 4 (a).  In the \textit{hot summer case}, the external ambient temperature is set to $35^{\circ}$C and the goal is to cool the cabin to the same target set point of $22^{\circ}$C by tracking the reference trajectory $T_{\mathrm{cab,ref}}$ depicted as a black line in Fig. 6 (a). In both scenarios the vehicle follows the WLTC driving cycle. The energy metric used for comparison is the electrical energy associated with the thermal subsystem only, computed by integrating the power absorbed by all the devices involved in the thermal management system.
State and input constraints are enforced over the prediction horizon. For
$x(k)=[T_{\mathrm{bat}}(k),\,SOC(k),\,\Delta T(k),\,T_{\mathrm{cab}}(k)]^{\top}$, the state bounds are
\begin{equation}
\begin{aligned}
-10^{\circ}\mathrm{C} &\le T_{\mathrm{bat}}(k{+}i|k)\le 36^{\circ}\mathrm{C},\\
-20^{\circ}\mathrm{C} &\le \Delta T(k{+}i|k)\le 20^{\circ}\mathrm{C},
\end{aligned}
\qquad i=1,\dots,p,
\label{eq:constr_states}
\end{equation}
while actuator saturation is modeled as
\begin{equation}
\begin{aligned}
0 &\le v_{\mathrm{comp}}(k{+}i|k)\le 8600\ \mathrm{rpm},\\
0.001 &\le \%\mathrm{valve}(k{+}i|k)\le 1,
\end{aligned}
\qquad i=1,\dots,p.
\label{eq:constr_inputs}
\end{equation}
%\vspace{-0.2cm}
The supervisory and piloting layers operate with a common sampling 
time of $T_s = 1$\,s, with prediction horizons of $p_{\mathrm{s}}=50$ 
and $p_{\mathrm{p}}=2$ steps respectively, reflecting the separation between long-horizon energy planning and short-horizon regulation.
The cost function weights are set to $Q=10^{8}$ and 
$R=[10^{15},\,10^{10}]$.
Figs.~\ref{fig:cold_states}--\ref{fig:cold_inputs} report the results obtained for the \textit{cold winter case}, and Figs.~\ref{fig:hot_states}--\ref{fig:hot_inputs} report the results obtained for the \textit{hot summer case}. In both conditions, the proposed control scheme (denoted as MPC) achieves improved tracking of $T_{\mathrm{cab,ref}}$ and maintains $T_{\mathrm{bat}}$ and $\Delta T$ within the admissible ranges. Compared to RB, the proposed strategy shapes compressor speed and valve routing more smoothly and anticipatively over the WLTC, reducing unnecessary actuation and the associated thermal energy demand, as clearly stated by comparison in terms of power consumption reported in Fig 4 (c) and Fig. 6 (c). The resulting energy savings with respect to the Stellantis RB strategy, computed by numerically integrating the power profile in Fig. 4 and 6., are summarized in Table~\ref{tab:energy-saving}: the proposed  control approach (denoted as MPC)  yields a reduction of 12.02\% in the cold winter case and 12.71\% in the hot summer case.

\begin{table}[t]
\centering
\caption{Thermal energy saving of MPC with respect to the Stellantis RB controller (WLTC).}
\label{tab:energy-saving}
\renewcommand{\arraystretch}{1.15}
\begin{tabular}{lcc}
\hline
\textbf{Scenario} & \textbf{MPC vs. RB} & \textbf{Ambient $T_{\mathrm{amb}}$} \\
 & \textbf{energy saving [\%]} & \textbf{[$^{\circ}$C]} \\
\hline
Cold case (heating)  & 12.02 & $-10$ \\
Hot case (cooling)   & 12.71 & $35$ \\
\hline
\end{tabular}
\end{table}

\section{Conclusion}
\label{sec:conclusion}

This paper presented a hierarchical NMPC strategy for integrated cabin HVAC and battery thermal management in BEVs. Neural-network-based prediction models for heating and cooling operation were identified using a CMO-based framework and embedded in a two-layer MPC architecture that combines long-horizon energy-oriented planning with short-horizon command tracking and corrective regulation. WLTC simulations under cold and hot ambient conditions show improved cabin temperature reference tracking and reduced thermal energy consumption, with energy savings of more than 12\% compared to an industrial rule-based benchmark, while also guaranteeing satisfaction of all operational constraints.
Future extensions include the investigation of robust MPC formulations to explicitly account for uncertainties in the 
previewed vehicle speed profile, as well as its inclusion as an additional decision variable within the optimization problem.
\vspace{-0.4cm}

\bibliographystyle{IEEEtran}
 \bibliography{biblio}

\end{document}